# DISCUSSION: CONDITIONAL GROWTH CHARTS

By Matias Salibian-Barrera[1] and Ruben H. Zamar[2]

*University of British Columbia*

The authors are to be congratulated for a very important contribution with many practical applications. Including covariates in the construction of growth charts will undoubtedly lead to more informative tools for practitioners in many disciplines.

Growth charts are widely used in practice to monitor the evolution of particular univariate measurements over time. In some situations, a patient's evolution may be better described by the joint behavior of these variables of interest. For example, instead of using two univariate growth charts to map the weight and height of children, the physician may prefer to simultaneously locate the child's measurements with respect to the joint population distribution of weight and height for children of the same age cohort. It is well known that biological variables are generally correlated and that these correlations may be important to determine multivariate boundaries for the "normal" range of the response vector. When measurements are correlated, univariate growth charts may be unable to represent interesting combined features of the variables of interest.

We can identify the following challenges when one tries to develop multivariate growth charts:

(a) choosing an appropriate definition of multivariate quantiles;

(b) modeling multivariate quantiles to include the subject's prior development history and other covariates; and

(c) devising visualization tools to display individual trajectories with respect to the reference populations.

Regarding item (a), a nice unified presentation of several definitions of multivariate quantiles along with an insightful account of desirable properties is given in [2]. A proper extension of Wei and He's model to the multivariate

---

Received December 2005.

[1]Supported by an NSERC Discovery grant.

[2]Supported by an NSERC Discovery grant.

setting [which would address (b) above] is of great interest but beyond the scope of this note. We will focus our discussion on item (c) for the simple case where the only covariate is time.

For simplicity of presentation, in what follows we will restrict our attention to the bivariate case and use quantiles based on Tukey's half-space depth [3]. For a random sample $\mathbf{x}_1, \ldots, \mathbf{x}_n$, depth-based multivariate quantiles can be obtained as follows:

(i) first, for each observation $\mathbf{x}_i$ compute its half-space depth $HD(\mathbf{x}_i)$, which is defined as the smallest fraction of observations included in a closed half-space with boundary line that passes through $\mathbf{x}_i$;

(ii) define the multivariate quantile of $\mathbf{x}_i$ as the sample quantile of $HD(\mathbf{x}_i)$ among all the $n$ half-space depths $HD(\mathbf{x}_j)$, $j = 1, \ldots, n$.

Figure 1 shows a bivariate dataset where "extreme" points (those with quantiles smaller than 0.05) are indicated with circles, and innermost points (corresponding to quantiles larger than 0.95) are shown with ×'s. Note that this definition is closely related to the following population definition of depth-based median-oriented multivariate quantiles [2]. The half-space depth of a

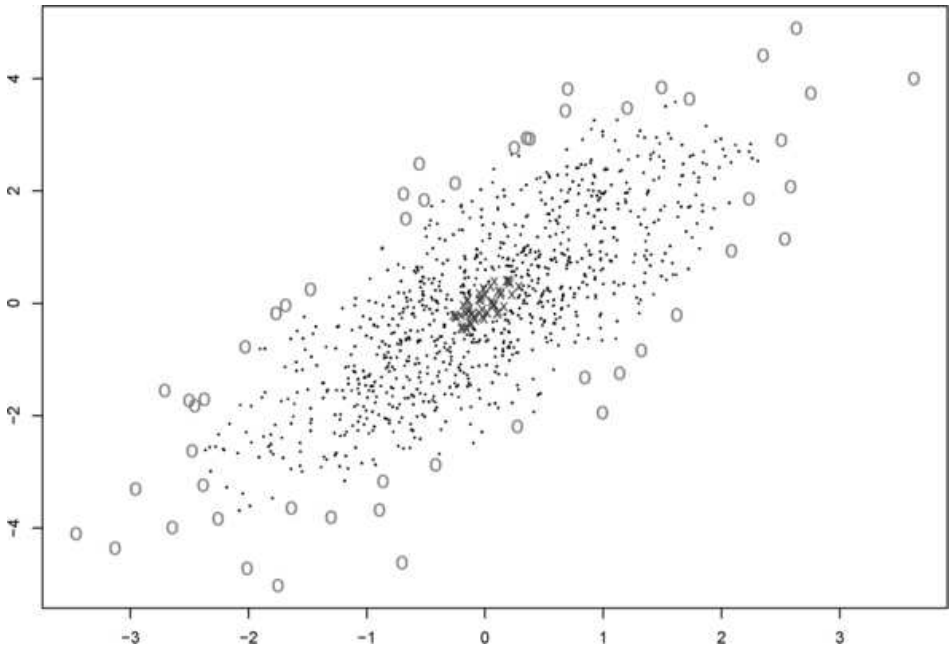

FIG. 1.  *A random sample of size $n = 1000$ from a bivariate normal distribution. Points with depth-based multivariate quantiles lower than 0.05 are shown with ∘'s, and those with quantiles higher than 0.95 are displayed with ×'s.*

vector $\mathbf{x}$ with respect to a reference probability measure $P$ is defined as

$$HD(\mathbf{x}, P) = \inf\{P(H) : H \text{ is a closed half-space and } \mathbf{x} \in H\}.$$

Following [2], the depth-based median-oriented multivariate $p$th quantile is the boundary of the set $I(\gamma_p, P)$, where

$$\gamma_p = \sup\{\gamma \geq 0 : P(I(\gamma, P)) \geq p\}$$

and

$$I(\gamma, P) = \{\mathbf{x} : HD(\mathbf{x}, P) \geq \gamma\}.$$

Multivariate quantiles are, of course, difficult to visualize when the dimension of the data is higher than 2 or 3. Since the main objective of growth charts is to locate the trajectory along time of a particular individual with respect to the corresponding different reference populations, we propose to use a series of univariate plots, one for each time point where the individual of interest has been observed.

More specifically, assume that we observe a vector of $p$ measurements $\mathbf{x}(t) \in \mathbb{R}^p$, at $k$ time points $t_1, \ldots, t_k$. Denote these observations by $\mathbf{x}(t_i)$, $i = 1, \ldots, k$. Furthermore, assume that we have $k$ reference populations from the distribution of the vector of measurements of interest $\mathbf{X}(t)$ at the same times $t_i$, $i = 1, \ldots, k$. For each time $t_i$, let $q_{t_i}$ be the corresponding multivariate quantile of the observed $\mathbf{x}(t_i)$ with respect to the corresponding reference populations. Let $\mathbf{a} \in \mathbb{R}^p$ with $\|\mathbf{a}\| = 1$ be any unit vector, and consider the projections $\mathbf{a}'\mathbf{x}(t_i)$, $i = 1, \ldots, k$. Denote by $\tilde{q}_{t_i}(\mathbf{a})$ the corresponding (univariate) quantiles of the projections $\mathbf{a}'\mathbf{x}(t_i)$ with respect to the projected reference populations. We propose to find the vector $\mathbf{a}_0$ for which the resulting $\tilde{q}_{t_i}(\mathbf{a}_0)$'s are closest to the multivariate quantiles $q_{t_i}$, $i = 1, \ldots, k$. Because in many applications the observations $\mathbf{x}(t_1), \ldots, \mathbf{x}(t_k)$ correspond to measurements taken over time on a particular patient, we will call this optimal vector $\mathbf{a}_0$ the "patient-specific direction." In other words, we define

$$(1.1) \qquad \mathbf{a}_0 = \underset{\|\mathbf{a}\| = 1}{\arg\min} \sum_{j=1}^{k} [q_{t_j} - \tilde{q}_{t_j}(\mathbf{a})]^2,$$

and use this direction to find univariate reference populations for which the corresponding projection of the measurements of the "patient of interest" are closest to the multivariate ones. Note that this patient-specific optimal direction (common for all time points) may provide some insight into which combination of the patient's measurements best describes the relative position of this individual in the multivariate reference populations. For example, if this patient's multivariate quantiles are becoming more "extreme" with time, the patient-specific optimal direction $\mathbf{a}_0$ may help the physician understand in which way this patient is deviating from the bulk of the reference populations.

We can now use different graphical tools to display the relative position of the individual of interest with respect to the projected reference populations. In the example below we use boxplots. Alternatively, one could use histograms or kernel density estimators.

EXAMPLE. We will illustrate the main ideas of our visualization proposal with a synthetic example. Although our proposal is aimed at situations involving several variables, in what follows we only consider bivariate observations to be able to plot the multivariate data and their projections.

We generated $k = 4$ reference samples of size 1000 from a bivariate normal distribution with constant correlation equal to 0.77 and variances 1 and 2.44. The mean vectors change over time from $(5,5)'$ to $(10.5,9)'$. In Figure 2 we display the four datasets along with the bivariate measurements corresponding to a "patient of interest," indicated with solid circles. The half-space depth bivariate quantiles for this patient are 0.79, 0.38, 0.07 and 0.02 at times $t_1$, $t_2$, $t_3$ and $t_4$, respectively. Clearly this patient is becoming atypical as time progresses.

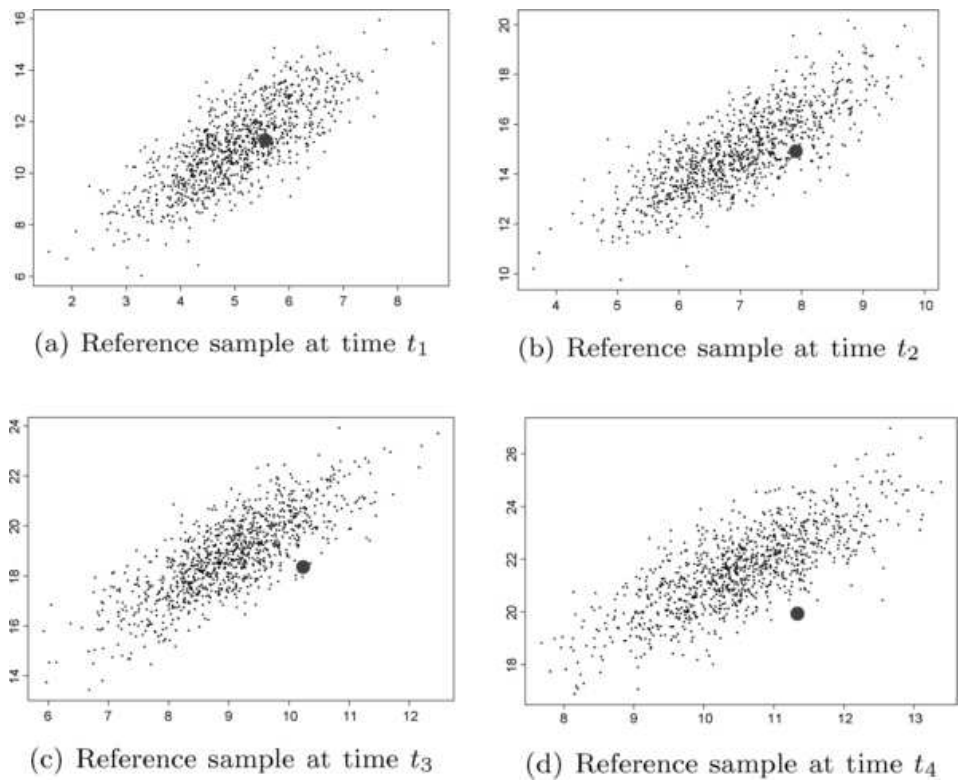

(a) Reference sample at time $t_1$

(b) Reference sample at time $t_2$

(c) Reference sample at time $t_3$

(d) Reference sample at time $t_4$

FIG. 2. *Reference samples of size* 1000, *at times* $t_1$, $t_2$, $t_3$ *and* $t_4$. *The measurements for a "patient of interest" are indicated with solid circles.*

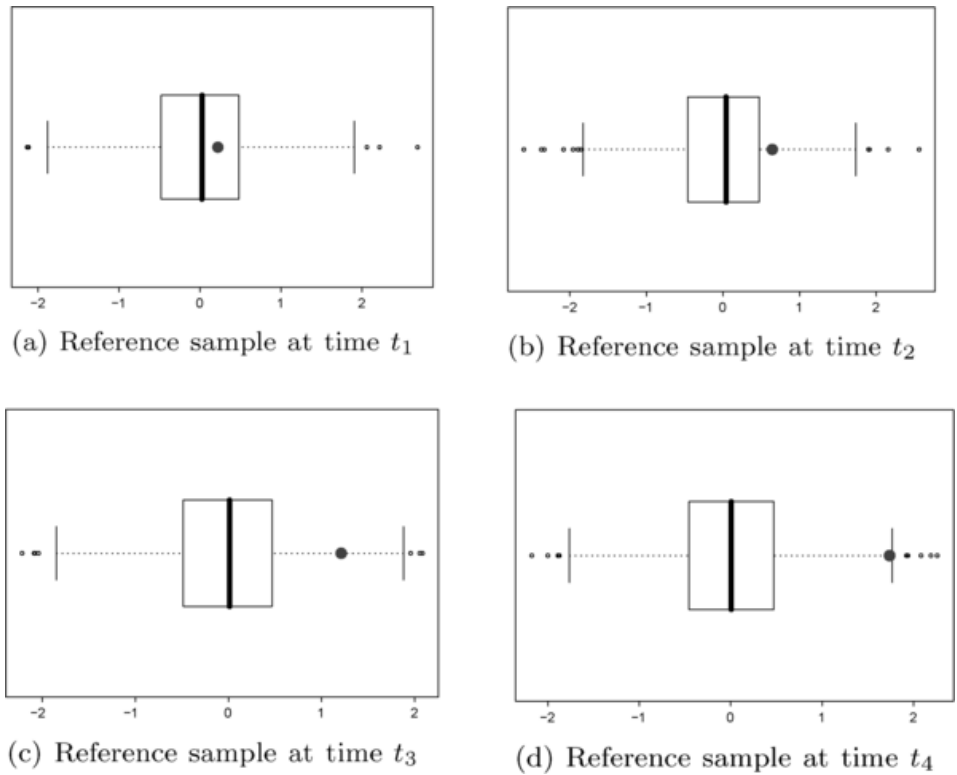

(a) Reference sample at time $t_1$      (b) Reference sample at time $t_2$

(c) Reference sample at time $t_3$      (d) Reference sample at time $t_4$

Fig. 3. *Relative position of the "patient of interest" (indicated with a solid circle) with respect to the reference samples projected on the patient-specific optimal direction.*

To find $\mathbf{a}_0$ in (1.1) we used a grid of 500 directions with equally spaced angles between 0 and $\pi$. Note that the multivariate quantiles $q_{t_j}$ only need to be computed once, regardless of the number of directions used in the numerical optimization. The bivariate half-space depths were computed using the AS 307 algorithm [1]. The optimal $\hat{\mathbf{a}}_0 = (0.72, -0.69)'$ and the univariate quantiles $\tilde{q}_{t_j}(\hat{\mathbf{a}}_0)$ were found to be 0.79, 0.38, 0.09 and 0.02, in close agreement with the multivariate ones. It is interesting to note that the "source" of this patient's increasing "unusualness" seems to be explained by a contrast between the first and second measurements. If these were height and weight, for example, we could conclude that this patient is becoming disproportionally tall for his/her weight.

Figure 3 shows the position of the "patient of interest" with respect to the reference samples projected on the patient-specific optimal direction $\hat{\mathbf{a}}_0$. The corresponding projections of the point of interest are indicated with a solid dot. We can clearly see that this point is noticeably deviating from the bulk of the population (crossing quantile boundaries) as time progresses. Furthermore, Figure 4 shows the relative position of the "patient of interest" using the two original variables. Note that, as it is to be expected, the univariate approach fails to detect the unusual behavior of this patient's measurements. This last figure illustrates the potential usefulness of multivariate growth charts compared with the analysis of several univariate growth charts.

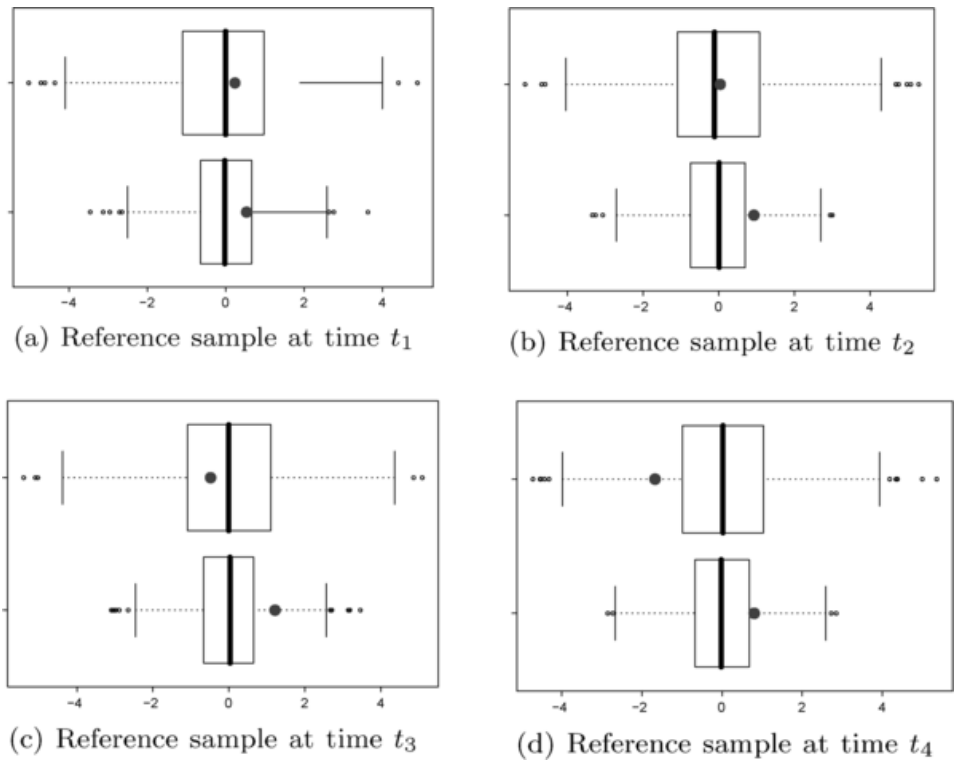

(a) Reference sample at time $t_1$

(b) Reference sample at time $t_2$

(c) Reference sample at time $t_3$

(d) Reference sample at time $t_4$

Fig. 4.   *Relative position of the "patient of interest" (indicated with a solid circle) with respect to the two univariate coordinates of the reference samples.*

## REFERENCES

[1] ROUSSEEUW, P. J. and RUTS, I. (1996). Algorithm AS 307: Bivariate location depth. *Appl. Statist.* **45** 516–526.

[2] SERFLING, R. (2002). Quantile functions for multivariate analysis: Approaches and applications. *Statist. Neerlandica* **56** 214–232. MR1916321

[3] TUKEY, J. W. (1975). Mathematics and the picturing of data. In *Proc. International Congress of Mathematicians* **2** (R. James, ed.) 523–531. Canadian Math. Congress, Montreal. MR0426989

DEPARTMENT OF STATISTICS
333-6356 AGRICULTURAL ROAD
UNIVERSITY OF BRITISH COLUMBIA
VANCOUVER, BRITISH COLUMBIA
CANADA V6T 1Z2
E-MAIL: ruben@stat.ubc.ca
        matias@stat.ubc.ca